\documentclass[10pt]{amsart}
\usepackage{ amssymb,latexsym, amscd,pb-diagram}
\usepackage{sseq}
\usepackage[all]{xy}
\usepackage[square, numbers]{natbib}
\usepackage{graphicx}
\usepackage{mathrsfs}
\usepackage{color}
\usepackage{amsmath}
\usepackage{multirow}

\usepackage{tikz}
\usepackage{dsfont}
\usetikzlibrary{matrix}

 \newtheorem{theorem}{Theorem}[section]

 \theoremstyle{definition}
 
 \theoremstyle{definition}

 \numberwithin{equation}{section}

\newcommand{\ben}{\begin{equation}}
\newcommand{\een}{\end{equation}}

\newcommand{\integer}{\ensuremath{{\mathbb Z}}}

\newcommand{\complex}{\ensuremath{{\mathbb C}}}

\newcommand{\Cx}{\ensuremath{{\mathbb C}^*}}

\newcommand{\field}{\ensuremath{{\mathbb F}}}

\newcommand{\DD}{{\mathcal D}}

\newcommand{\ZZ}{{\mathcal Z}}

\newcommand{\VV}{{\mathcal V}}

\newcommand{\CC}{\mathcal{C}}

\newcommand{\MM}{\mathcal{M}}

\newcommand{\Hom}{\mathrm{Hom}}

\newcommand{\Map}{\ensuremath{{\mathrm{Map}}}}

\newcommand{\Z}{{\mathbb Z}}

\newcounter{commentcounter}

\begin{document}

\title[Classification of Pointed Fusion Categories of dimension 8]{Classification of Pointed Fusion Categories of dimension 8 up to weak Morita Equivalence}

\thanks{The first author acknowledges the support 
of COLCIENCIAS through grant number FP44842-087-2017 of the Convocatoria Nacional J\'ovenes Investigadores e Innovadores No 761 de 2016. The second author acknowledges the financial support of the Max Planck Institute of Mathematics in Bonn, Germany.}
%\author{C\'esar Galindo}
%\address{Departamento de Matem\'{a}ticas, Universidad de los Andes,
%Carrera 1 N. 18A - 10, Bogot\'a, Colombia.  }
%\email{
%cn.galindo1116@uniandes.edu.co, cesarneyit@gmail.com}
\author{\'Alvaro Mu\~noz}
\address{Departamento de Matem\'{a}ticas y Estad\'istica, Universidad del Norte, Km.5 V\'ia Antigua a Puerto Colombia, 
Barranquilla, Colombia.}
\email{munozfa@uninorte.edu.co}
\author{Bernardo Uribe}
\email{bjongbloed@uninorte.edu.co}
%\date{June 20th  2013}
\begin{abstract}
In this paper we give a complete classification of pointed fusion categories over $\complex$ of global dimension 8.
We first classify the equivalence classes of pointed fusion categories of dimension 8, and then we 
proceed to determine which of these equivalence classes have equivalent categories of modules. 
This classificaction permits to classify the equivalence classes of braded tensor equivalences
of twisted Drinfeld doubles of finite groups of order 8.
\end{abstract}

\maketitle

\section*{Introduction}
A fusion category \cite{ENO} is a rigid semisimple $\complex$-linear tensor category with only finitely many isomorphism classes of simple objects, such that the endomorphisms of the unit object is $\complex$. It is moreover pointed if 
all its  simple objects are invertible. Any pointed fusion category is equivalent to a fusion category $Vect(H, \eta)$ of $H$-graded
complex vector spaces with $H$ a finite group, together with an associativity constraint defined by a cocycle $\eta \in Z^3(H, \complex^*)$.
The skeletal tensor category $\VV(H, \eta)$ associated to $Vect(H, \eta)$ is a 2-group \cite{BaezLauda} 
with only
one object for each isomorphism class of simple objects in $Vect(H, \eta)$, whose
morphisms are only automorphisms and whose associator is defined by $\eta$. 
For any module category $\MM$ over a fusion category $\CC$ we may define the dual fusion category $\CC_\MM^*:=Fun_\CC(\MM,\MM)$, and we say that two fusion categories $\CC$ and $\DD$ are weakly
Morita equivalent if there exists an idecomposable module category $\MM$ over $\CC$ such that $\CC_\MM^*$ and $\DD$
are tensor equivalent \cite[Def. 4.2]{MugerI}. 

In \cite{Naidu} there were given necessary and sufficient conditions in terms of cocycles for two 
pointed fusion categories $Vect(H, \eta)$ and $Vect(\widehat{H}, \widehat{\eta})$ to be weakly Morita equivalent,
and in \cite{Uribe} a choice of appropriate coordinates permitted the second author to give a precise description
of the groups $H, \widehat{H}$ and the cocycles $\eta, \widehat{\eta}$ for the pointed fusion categories to be weakly
Morita equivalent.
In this paper we follow the description done by the second author in \cite{Uribe} in order to classify 
the Morita equivalence classes of pointed fusion categories of global dimension 8

This work will be divided in two chapters. In the first chapter we will setup explicit basis for $H^3(H, \complex^*)$ 
 and we will calculate the space of orbits $H^3(H, \complex^*)/Aut(H)$ for each of the five groups of order 8; this
 will determine the equivalence classes of pointed fusion categories of global dimension 8. In the second 
 chapter we will recall the classification theorem of pointed fusion categories \cite[Thm. 3.9]{Uribe} and we will
 use the Lyndon-Hochschild-Serre spectral sequence associated to group extensions to explicitly find the Morita equivalence classes
 of pointed fusion categories of global dimension 8. We will finish the paper with an application to the
 classification of braided tensor equivalences of twisted Drinfeld doubles of groups of order 8.
 
 Since this work is an application of the results obtained by the second author in \cite{Uribe}, we 
 will use the notation and the constructions done there. A more detailed version of the results presented in this paper appear in the Master thesis of the first author \cite{Munoz}.

\section{Equivalence classes of pointed fusion categories of global dimension 8}

The pointed fusion categories $Vect(H, \eta)$ and $Vect(\widehat{H}, \widehat{\eta})$ are equivalent 
if and only if there is and isomorphism of groups $\phi: H \stackrel{\cong}{\to} \widehat{H}$ such that
$[\phi^* \widehat{\eta}]=[\eta]$ in $H^3(H, \complex^*)$. Therefore the equivalence classes
of pointed fusion categories of global dimension 8 is isomorphic to the union of the spaces of orbits $H^3(H, \complex^*)/Aut(H)$
where $H$ runs over the groups $\integer_2^3, \integer_4 \times \integer_2, \integer_8,D_8$ and $Q_8$.
By the short exact sequence of coefficients $0 \to \integer \to \complex \to \complex^* \to 0$ we now that for any finite group $H^3(H, \complex^*) \cong  H^4(H, \integer)$;
 in what follows we will calculate $H^4(H, \integer)/Aut(H)$.
\subsection{$(\integer/2)^3$}
We know that
$$H^4((\integer/2)^3, \integer) = ker \left( Sq^1 : H^4((\integer/2)^3, \field_2) \to H^5((\integer/2)^3, \field_2)
\right)$$
where $\field_2$ is the field of 2 elements and $Sq^1$ is the first Steenrod square operation in cohomology with $\field_2$ coefficients. Letting $H^*((\integer/2)^3, \field_2) = \field_2[x,y,z]$ we obtain
$$H^4((\integer/2)^3, \integer) \cong\integer/2 \langle x^4,y^4,z^4,x^2y^2,x^2z^2,y^2z^2,x^2yz+xy^2z+xyz^2 \rangle
\cong (\integer/2)^7.$$

The group of automorphisms of $(\integer/2)^3$ is isomorphic to $GL(3, \field_2)$ and the 10 orbits of
$H^4((\integer/2)^3, \integer) / Aut((\integer/2)^3)$ are:
\begingroup
\everymath{\scriptstyle}
\scriptsize
\begin{align*}orb(0) =\{0\}\end{align*}
\begin{align*}orb(x^4) =\{ x^{4}, y^{4}, z^{4}, x^4+y^{4}, x^{4}+z^{4}, y^{4}+z^{4},x^{4}+y^{4}+z^{4}\}\end{align*}
\begin{align*}
orb(&x^2y^2) = \\
\{
& x^{2}y^{2}, x^{2}y^{2}+y^{4}, x^{2}z^{2}+x^{4}, x^{2}y^{2}+y^{2}z^{2}, x^{4}+x^{2}y^{2}+x^{2}z^{2}+y^{2}z^{2}, x^{4}+x^{2}y^{2}+x^{2}z^{2},\\
&x^{2}z^{2}, x^{2}y^{2}+x^{4}, y^{2}z^{2}+y^{4}, x^{2}y^{2}+x^{2}z^{2},y^{4}+x^{2}y^{2}+x^{2}z^{2}+y^{2}z^{2}, y^{4}+x^{2}y^{2}+y^{2}z^{2},\\
&y^{2}z^{2},x^{2}z^{2}+z^{4}, y^{2}z^{2}+z^{4}, x^{2}z^{2}+y^{2}z^{2}, z^{4}+x^{2}y^{2}+x^{2}z^{2}+y^{2}z^{2}, z^{4}+x^{2}z^{2}+y^{2}z^{2},\\
&x^{4}+y^{4}+x^{2}z^{2}+y^{2}z^{2},
z^{4}+y^{4}+x^{2}y^{2}+x^{2}z^{2},
x^{4}+z^{4}+x^{2}y^{2}+y^{2}z^{2} \}
\end{align*}
\begin{align*}
orb(& x^{4}+y^{2}z^{2}) = \\
\{
& x^{4}+y^{2}z^{2}, x^{4}+y^{4}+x^{2}y^{2}+x^{2}z^{2}, x^{4}+z^{4}+x^{2}z^{2}+y^{2}z^{2}, x^{4}+y^{4}+x^{2}z^{2}, x^{4}+z^{4}+x^{2}y^{2},  \\
&y^{4}+x^{2}z^{2}, x^{4}+y^{4}+x^{2}y^{2}+y^{2}z^{2}, y^{4}+z^{4}+x^{2}y^{2}+y^{2}z^{2}, x^{4}+y^{4}+y^{2}z^{2}, y^{4}+z^{4}+x^{2}y^{2},  \\
&z^{4}+x^{2}y^{2}, x^{4}+z^{4}+x^{2}y^{2}+x^{2}z^{2}, y^{4}+z^{4}+x^{2}z^{2}+y^{2}z^{2}, x^{4}+z^{4}+z^{2}y^{2}, y^{4}+z^{4}+x^{2}z^{2}, \\
&x^{4}+y^{4}+z^{4}+x^{2}y^{2}, x^{4}+x^{2}y^{2}+y^{2}z^{2}, z^{4}+x^{2}y^{2}+y^{2}z^{2},
x^{4}+y^{4}+x^{2}y^{2}+x^{2}z^{2}+y^{2}z^{2},\\
&x^{4}+y^{4}+z^{4}+x^{2}z^{2}, x^{4}+x^{2}z^{2}+y^{2}z^{2}, y^{4}+x^{2}y^{2}+x^{2}z^{2}, 
x^{4}+z^{4}+x^{2}y^{2}+x^{2}z^{2}+y^{2}z^{2}, \\
&x^{4}+y^{4}+z^{4}+y^{2}z^{2},  z^{4}+x^{2}y^{2}+x^{2}z^{2} , y^{4}+x^{2}z^{2}+y^{2}z^{2},
y^{4}+z^{4}+x^{2}y^{2}+x^{2}z^{2}+y^{2}z^{2}, \\
&x^{2}y^{2}+x^{2}z^{2}+y^{2}z^{2}
  \}
\end{align*}
\begin{align*}
orb(& x^{4}+y^{4}+x^{2}y^{2}) = \\
\{
&x^{4}+y^{4}+x^{2}y^{2}, x^{4}+y^{4}+z^{4}+x^{2}y^{2}+x^{2}z^{2}, 
x^{4}+z^{4}+x^{2}z^{2}, x^{4}+y^{4}+z^{4}+x^{2}y^{2}+y^{2}z^{2}, \\
&y^{4}+z^{4}+y^{2}z^{2}, x^{4}+y^{4}+z^{4}+x^{2}z^{2}+y^{2}z^{2},
x^{4}+y^{4}+z^{4}+x^{2}y^{2}+x^{2}z^{2}+y^{2}z^{2}
  \}
\end{align*}
\begin{align*}
orb(& x^{4}+x^{2}yz+xy^{2}z+xyz^{2}) = \\
\{
&x^{4}+x^{2}yz+xy^{2}z+xyz^{2},  x^{4}+y^{4}+z^{4}+x^{2}yz+xy^{2}z+xyz^{2}, x^{4}+z^{4}+x^{2}yz+xy^{2}z+xyz^{2}+y^{2}z^{2},\\
&y^{4}+x^{2}yz+xy^{2}z+xyz^{2}, x^{4}+y^{4}+x^{2}yz+xy^{2}z+xyz^{2}+y^{2}z^{2}, y^{4}+x^{2}yz+xy^{2}z+xyz^{2}+y^{2}z^{2},\\
&z^{4}+x^{2}yz+xy^{2}z+xyz^{2}, z^{4}+x^{2}yz+xy^{2}z+xyz^{2}+y^{2}z^{2}, y^{4}+x^{2}yz+xy^{2}z+xyz^{2}+x^{2}y^{2},\\
&x^{4}+z^{4}+x^{2}yz+xy^{2}z+xyz^{2}+x^{2}y^{2}, y^{4}+z^{4}+x^{2}yz+xy^{2}z+xyz^{2}+x^{2}z^{2},\\
&y^{4}+z^{4}+x^{2}yz+xy^{2}z+xyz^{2}+x^{2}y^{2}, x^{4}+y^{4}+x^{2}yz+xy^{2}z+xyz^{2}+x^{2}z^{2},\\
&x^{4}+x^{2}yz+xy^{2}z+xyz^{2}+x^{2}y^{2}, z^{4}+x^{2}yz+xy^{2}z+xyz^{2}+x^{2}z^{2}, x^{4}+z^{4}+x^{2}yz+xy^{2}z+xyz^{2}+x^{2}y^{2}+y^{2}z^{2}\\
&x^{4}+x^{2}yz+xy^{2}z+xyz^{2}+x^{2}z^{2}, z^{4}+x^{2}yz+xy^{2}z+xyz^{2}+x^{2}y^{2}+x^{2}z^{2}\\
&y^{4}+z^{4}+x^{2}yz+xy^{2}z+xyz^{2}+x^{2}y^{2}+x^{2}z^{2}, y^{4}+x^{2}yz+xy^{2}z+xyz^{2}+x^{2}y^{2}+x^{2}z^{2}\\
&x^{2}yz+xy^{2}z+xyz^{2}+x^{2}y^{2}+x^{2}z^{2},  x^{4}+x^{2}yz+xy^{2}z+xyz^{2}+x^{2}y^{2}+y^{2}z^{2}\\
&x^{2}yz+xy^{2}z+xyz^{2}+x^{2}y^{2}+y^{2}z^{2}, y^{4}+x^{2}yz+xy^{2}z+xyz^{2}+x^{2}z^{2}+y^{2}z^{2}, \\
&z^{4}+x^{2}yz+xy^{2}z+xyz^{2}+x^{2}y^{2}+y^{2}z^{2}, x^{4}+y^{4}+xyz+xy^{2}z+xyz^{2}+x^{2}z^{2}+y^{2}z^{2}, \\
&xyz+xy^{2}z+xyz^{2}+x^{2}z^{2}+y^{2}z^{2}, x^{4}+xyz+xy^{2}z+xyz^{2}+x^{2}z^{2}+y^{2}z^{2}  \}
\end{align*}

\begin{align*}
orb(&x^{4}+x^{2}yz+xy^{2}z+xyz^{2}+y^{2}z^{2}) = \\
\{
&x^{4}+x^{2}yz+xy^{2}z+xyz^{2}+y^{2}z^{2}, x^{4}+z^{4}+x^2yz+xy^{2}z+xyz^{2}+x^{2}z^{2}+y^{2}z^{2},\\
&y^{4}+x^{2}yz+xy^{2}z+xyz^{2}+x^{2}z^{2}, y^{4}+z^{4}+x^2yz+xy^{2}z+xyz^{2}+x^{2}z^{2}+y^{2}z^{2},\\
&z^{4}+x^{2}yz+xy^{2}z+xyz^{2}+x^{2}y^{2}, x^{4}+z^{4}+x^2yz+xy^{2}z+xyz^{2}+x^{2}y^{2}+x^{2}z^{2},\\
&x^{4}+y^{4}+x^2yz+xy^{2}z+xyz^{2}+x^{2}y^{2}+x^{2}z^{2}, x^{4}+y^{4}+x^{2}yz+xy^{2}z+xyz^{2},\\
&x^{4}+y^{4}+x^2yz+xy^{2}z+xyz^{2}+x^{2}y^{2}+y^{2}z^{2}, x^{4}+z^{4}+x^{2}yz+xy^{2}z+xyz^{2},\\
&y^{4}+z^{4}+x^2yz+xy^{2}z+xyz^{2}+x^{2}y^{2}+y^{2}z^{2}, y^{4}+z^{4}+x^{2}yz+xy^{2}z+xyz^{2},\\
&x^{4}+z^{4}+x^{2}yz+xy^{2}z+xyz^{2}+x^{2}z^{2}, x^{4}+y^{4}+z^{4}+x^{2}yz+xy^{2}z+xyz^{2}+x^{2}y^{2},\\
&x^{4}+y^{4}+x^{2}yz+xy^{2}z+xyz^{2}+x^{2}y^{2}, x^{4}+y^{4}+z^{4}+x^{2}yz+xy^{2}z+xyz^{2}+x^{2}z^{2},\\
&y^{4}+z^{4}+x^{2}yz+xy^{2}z+xyz^{2}+y^{2}z^{2}, x^{4}+y^{4}+z^{4}+x^{2}yz+xy^{2}z+xyz^{2}+y^{2}z^{2},\\
&x^{4}+y^{4}+z^{4}+x^{2}yz+xy^{2}z+xyz^{2}+x^{2}y^{2}+x^{2}z^{2},\\
&x^{4}+y^{4}+z^{4}+x^{2}yz+xy^{2}z+xyz^{2}+x^{2}y^{2}+y^{2}z^{2},\\
&x^{4}+y^{4}+z^{4}+x^{2}yz+xy^{2}z+xyz^{2}+x^{2}z^{2}+y^{2}z^{2}
  \}
\end{align*}
\begin{align*}
orb(& x^{2}yz+xy^{2}z+xyz^{2}+x^{2}y^{2}+x^{2}z^{2}+y^{2}z^{2}) = \\
\{
 & x^{2}yz+xy^{2}z+xyz^{2}+x^{2}y^{2}+x^{2}z^{2}+y^{2}z^{2},\\
& x^{4}+z^{4}+x^{2}yz+xy^{2}z+xyz^{2}+x^{2}y^{2}+x^{2}z^{2}+y^{2}z^{2}, 
x^{4}+x^{2}yz+xy^{2}z+xyz^{2}+x^{2}y^{2}+x^{2}z^{2}+y^{2}z^{2},\\
&y^{4}+z^{4}+x^{2}yz+xy^{2}z+xyz^{2}+x^{2}y^{2}+x^{2}z^{2}+y^{2}z^{2},
z^{4}+x^{2}yz+xy^{2}z+xyz^{2}+x^{2}y^{2}+x^{2}z^{2}+y^{2}z^{2},\\
&x^{4}+y^{4}+x^{2}yz+xy^{2}z+xyz^{2}+x^{2}y^{2}+x^{2}z^{2}+y^{2}z^{2},
y^{4}+x^{2}yz+xy^{2}z+xyz^{2}+x^{2}y^{2}+x^{2}z^{2}+y^{2}z^{2}
  \}
\end{align*}
\begin{align*}
orb(&x^{2}yz+ xy^{2}z+xyz^{2}) = \\
\{ &x^{2}yz+xy^{2}z+xyz^{2},\\
&x^{2}yz+xy^{2}z+xyz^{2}+y^{2}z^{2}, x^{4}+x^{2}yz+xy^{2}z+xyz^{2}+x^{2}y^{2}+x^{2}z^{2},\\
&x^{2}yz+xy^{2}z+xyz^{2}+x^{2}y^{2}, y^{4}+x^{2}yz+xy^{2}z+xyz^{2}+x^{2}y^{2}+y^{2}z^{2},\\
&x^{2}yz+xy^{2}z+xyz^{2}+x^{2}z^{2}, z^{4}+x^{2}yz+xy^{2}z+xyz^{2}+x^{2}z^{2}+y^{2}z^{2}
 \}
\end{align*}
\begin{align*}
orb(x^{4}+ y^{4}+z^{4}+x^{2}yz+ & xy^{2}z+xyz^{2}+x^{2}y^{2}+x^{2}z^{2}+y^{2}z^{2})\\
=& \{x^{4}+ y^{4}+z^{4}+x^{2}yz+xy^{2}z+xyz^{2}+ x^{2}y^{2}+x^{2}z^{2}+y^{2}z^{2} \}
\end{align*}
\endgroup

\subsection{$\integer/4 \times \integer/2$} \label{subsection Z4xZ2} By Kunneth's theorem we know that
$$H^4(\integer/4 \times \integer/2, \integer) \cong \integer \langle v^2, uv, u^2 \rangle /(4v^2,2uv, 2u^2) \cong
\integer/4 \oplus \integer/2 \oplus \integer/2$$
where $H^*(\integer/4,\integer)=\integer[v]/(4v)$ and $H^*(\integer/2,\integer)=\integer[u]/(4u)$. 
The group of automorphisms of $\integer/4 \times \integer/2$ is the dihedral group $D_8$ generated
by the automorphisms
\begin{align*}
\rho:\integer/4 \times \integer/2 &\longrightarrow\integer/4 \times \integer/2&\sigma:\integer/4 \times \integer/2&\longrightarrow \integer/4 \times \integer/2\\
(1,0)&\longrightarrow (1,1)&(1,0)&\longrightarrow (1,1)\\
(0,1)&\longrightarrow (2,1)&(0,1)&\longrightarrow (0,1)
\end{align*}
with $\sigma \rho \sigma= \rho^{-1}$. The induced action in $H^2(\integer/4 \times \integer/2,\integer)$
is given by the equations
\begin{align*}
\rho^*u=u+v, \ \ \
\rho^*v=2u+v, \ \ \
\sigma^*u=u, \ \ \
\sigma^*v=2u+v
\end{align*}
and therefore the 9 orbits are
\begin{align*}
H^4&(\integer/4 \times \integer/2,\integer)/Aut(\integer/4 \times \integer/2)=\\
\lbrace
&\lbrace u^2, u^2+v^2\rbrace,
\lbrace v^2\rbrace,
\lbrace uv, 2u^2+uv+v^2, 2u^2+uv, uv+v^2\rbrace,
\lbrace u^2+uv, 3u^2+uv\rbrace, \\
&\lbrace 3u^2, 3u^2+v^2\rbrace, 
\lbrace u^2+uv+v^2, 3u^2+uv+v^2\rbrace,
\lbrace 2u^2\rbrace,
\lbrace 2u^2+v^2\rbrace,\lbrace 0\rbrace
\rbrace.\end{align*}

\subsection{$\integer/8$} The multiplicative group $\integer/8^*$ of units of $\integer/8$
is isomorphic to the automorphism group $Aut(\integer/8)$  and acts on $H^*(\integer/8, \integer) = \integer[s]/(8s)$
by multiplication on the generator $s$. Since all the units of $\integer/8$ square to 1, the action
of the units on $H^4(\integer/8, \integer)=\integer \langle s^2\rangle /(8s^2)$ is trivial. Therefore
$$H^4(\integer/8, \integer)/Aut(\integer/8)= \integer \langle s^2\rangle /(8s^2).$$

\subsection{$D_8$ } \label{subsection dihedral} The dihedral group $D_8= \langle a, b | a^4=b^2=1, b a b = a^{-1} \rangle$
has for automorphisms a group isomorphic to $D_8$ generated by the automorphisms
$\phi, \theta \in Aut(D_8)$ with $\phi(a)=a, \phi(b)=ba$, $\theta(a)=a^{-1}$, $\theta(b)=b$,
satisfying $\theta \phi \theta=\phi^{-1}$. The inner automorphisms are generated by $ad_b=\theta$ and
$ad_a=\phi^2$, 
and since the inner automorphisms act trivially on the cohomology
of the group, we only need to calculate the orbits that the action of $\phi^*$ induce on $H^4(D_8, \integer)$.

Consider the short exact sequence of groups
$$\xymatrix{\langle a^{2}\rangle=\mathbb{Z}/{2}\ar[r]&D_{8}\ar[r]^(0.2){\pi}& \mathbb{Z}/{2}\times\mathbb{Z}/{2}\cong {\left\langle a\right\rangle}/_{(a^{2})}\times\left\langle \tau\right\rangle }$$
with $\pi(a)=(1,0)$ and $\pi(b)=(0,1)$, and define 
\begin{align*}
\widehat{\phi}:\mathbb{Z}/{2}\times\mathbb{Z}/{2}&\longrightarrow \mathbb{Z}/{2}\times\mathbb{Z}/{2}\\
(1,0)&\longrightarrow(1,0)\\
(0,1)&\longrightarrow(1,1).
\end{align*}
Note that the isomorphism $\widehat{\phi}$ fits into the diagram
$$\xymatrix{{\mathbb{Z}/{2}\ar[r]\ar@{=}[d]} & D_{8}\ar[d]^{\phi}\ar[r] & \mathbb{Z}/{2}\times\mathbb{Z}/{2}\ar[d]_{\widehat{\phi}}\\
\mathbb{Z}/{2}\ar[r] & D_{8}\ar[r] & \mathbb{Z}/{2}\times\mathbb{Z}/{2}}$$
since $\phi$ fixes the center. We will fix a base for $H^3(D_8, \complex^*)$ using
the Lyndon-Hochschild-Serre spectral sequence associated to the short exact sequence,
and since $\phi$ preserve the center, we can deduce
the action of $\phi^*$ by looking at its action on the spectral sequence.

Let $H^*(\integer/2 \times \integer/2, \field_2)= \field_2[x,y]$ and note that the
$k$-invariant of the extension of $D_8$ is $xy+x^2$ since on the restrition $D_8|_{\integer/2 \times \{0\}}$
we get the extension $0 \to \integer/2 \to \integer/4 \to \integer/2 \to 0$ whose $k$-invariant is $x^2$.
The action of $\widehat{\phi}^*$ on $\field_2[x,y]$ is given by the equations 
$\widehat{\phi}^*x = x+y$ and $\widehat{\phi}^*y =y$.

The second page of the Lyndon-Hochschild-Serre spectral sequence 
$$E_2^{p,q}=H^p(\integer/2 \times \integer/2, H^q(\integer/2, \complex^*))$$
has for relevant terms
\renewcommand{\arraystretch}{2}
$$\newcommand*{\tempb}{\multicolumn{1}{|c}{}}
\begin{array}{cccccccccc}
3 & \tempb & \langle z^4\rangle\\
2 & \tempb & 0& 0\\
1 & \tempb &\langle z^2\rangle&\langle zx, zy\rangle&\langle zx^2,zy^2,zxy\rangle\\
0 & \tempb & \mathbb{C}^*&\langle x^2, y^2\rangle & \langle x^2y+xy^2\rangle& \langle x^4,y^4, x^2y^2\rangle& \langle x^4y+x^3y^2, x^2y^3+xy^4\rangle  \\ \cline{2-7}
 & & 0& ~~~~1~~~~ & ~~~~2~~~~ & ~~~~3~~~~ & 4
\end{array}
$$
where we have taken in the fiber $H^*(\integer/2,\field_2)=\field_2[z]$, and on $E_2^{*,0}$
and $E_2^{0,*}$ we have taken the elements in $H^*(\integer/2 \times \integer/2, \integer)$
and $H^*(\integer/2, \integer)$ respectively wich are trivial under the action of $Sq^1$.
Note that $\phi^*$ leaves $z$ fixed and on $x,y$ acts through the induced action of $\widehat{\phi}$.

The element $z^4$ survives the spectral sequence and the second differential $d_2: H^p(\integer/2 \times \integer/2,
\field_2) \to H^{p+2}(\integer/2 \times \integer/2, \complex^*)$
is determined by the second differential $\overline{d}_2$ on the LHS spectral sequence with coefficients in $\field_2$.
In this case $\overline{d}_2(z)=xy +x^2$ is the $k$-invariant and the second differential satisfies the equation
$$d_2(z p(x,y))= Sq^1\left((xy+x^2)p(x,y)\right)$$
where $p(x,y)$ is any polynomial in $x$ and $y$.
The relevant terms of the third page of the spectral sequence are
\renewcommand{\arraystretch}{2}
$$
\newcommand*{\tempb}{\multicolumn{1}{|c}{}}
\begin{array}{cccccccccc}
3 & \tempb &\mathbb\langle z^4\rangle\\
2 & \tempb & 0& 0\\
1 & \tempb & 0 & \langle zy\rangle & \langle z(xy+x^2)\rangle \\
0 & \tempb & \mathbb{C}^*& \langle x^2,y^2\rangle & 0 &\dfrac{\langle x^4,x^2y^2,y^4\rangle}{\langle x^2y^2+x^4\rangle}& 0 \\ \cline{2-7}
& & 0~~~& 1~~~ & 2~~~ & 3~~~ & 4
\end{array}
$$
and they all survive to the page at infinity.

The cohomology of the dihedral group is $H^3(D_8,\complex^*) = \integer/4 \oplus \integer/2 \oplus \integer/2$
and we may denote the generators  
\begin{align} \label{generators D8}
H^3(D_8,\complex^*) = \langle \gamma \rangle \oplus \langle
\alpha \rangle \oplus \langle \beta \rangle
\end{align} with $\beta= x^4$, $\alpha = y^4$, $E_3^{3,0}\cong \langle \gamma \rangle/(2\gamma)$ and $E_3^{2,1}\cong \langle 2\gamma \rangle / (4 \gamma)$. This choice of base
agrees with the fact that when we restrict the spectral sequence to the subgroup $\langle a \rangle$ which fits in
the short exact sequence $0 \to \integer/2 \to \integer/4 \to \integer/2 \to 0$, the relevant surviving terms
 of the spectral sequence are only $E_2^{3,0}$ and $E_3^{2,1}$ (it is the same as making $y=0$).

Now, since $\phi^*$ leaves $z$ and $y$ fixed, and $\phi^*x=x+y$, we obtain that the 12 orbits of the action
of $\phi^*$ are the following:
 \begin{align*}
H^3(D_8,\complex^*) &/Aut(D_8)=\lbrace \{0\},\{\gamma\},\{2\gamma\},\{3\gamma\},\{\alpha\},\{\beta,\alpha+\beta\},\lbrace\gamma +\beta,\gamma+\alpha+\beta\rbrace,\\
&\lbrace 2\gamma +\beta,2\gamma+\alpha+\beta\rbrace,\lbrace 3\gamma +\beta,3\gamma+\alpha+\beta \rbrace,
\lbrace\gamma+\alpha\rbrace,
\lbrace 2\gamma+\alpha\rbrace, \lbrace 3\gamma+\alpha\rbrace \rbrace 
\end{align*}

\subsection{$Q_8$} The group of quaternions $Q_8 =\{\pm1, \pm i, \pm j, \pm k \}$ is a subgroup
of $SU(2)$ and therefore has periodic cohomology with $H^4(Q_8, \integer)= \integer/8$. The group
of automorphisms $Aut(Q_8)$ is isomorphic to the symetric group $\mathfrak{S}_4$.

The resolution
$$\cdots\longrightarrow\Z Q_8\oplus\Z Q_8\stackrel{\delta_1}{\longrightarrow}\Z Q_8\stackrel{\delta_4}{\longrightarrow}\Z Q_8\stackrel{\delta_3}{\longrightarrow}\Z Q_8\oplus\Z Q_8\stackrel{\delta_2}{\longrightarrow}\Z Q_8\oplus\Z Q_8\stackrel{\delta_1}{\longrightarrow}\Z Q_8\stackrel{\epsilon}{\longrightarrow}\Z
 $$
with differentials $\delta_i$ defined by the equations 
\begin{align*}
&\delta_1(a_1,a_2)=a_1(i-1)+a_2(j-1)\\
&\delta_2(a_1,a_2)=(a_1(i+1)+a_2(ij+1),-a_1(j+1)+a_2(i-1))\\
&\delta_3(a)=(a(i-1),-a(ij-1))\\
&\delta_4(a)=a\underset{\tau\in Q_8}{\sum }\tau
\end{align*}
provides a periodic resolution of $\integer$ by free $\integer Q_8$-modules. Applying the functor
 $Hom_{\integer Q_8}(\cdot, \integer)$ we obtain a cochain complex that calculates the cohomology of $Q_8$ and
 on degree 4 we obtain $Hom_{\integer Q_8}(\integer Q_8, \Z) $ which is invariant under the action of the automorphis
 group $Aut(Q_8)$. Therefore $Aut(Q_8)$ acts trivially on $H^4(Q_8, \integer)$ and we obtain 8 orbits
 $$H^4(Q_8, \integer)/Aut(Q_8) = H^4(Q_8, \integer)=\integer \langle t \rangle /(8t).$$
 
 \vspace{0.5cm}
 
 We conclude that are 10+9+8+12+8=47 equivalence classes of pointed fusion categories of global dimension 8.

\section{Morita equivalence classes of pointed fusion categories of global dimension 8}

In this section we will explicitly determine which pointed fusion categories of global dimension 8 have equivalent
categories of modules categories. We will use the notation and the results of \cite{Uribe} and we will recall
the classification theorem \cite[Thm. 3.9]{Uribe}.

An skeletal indecomposable module category $\MM=(A \backslash H,\mu)$ of $\CC=\VV(H, \eta)$ is determined by a transitive $H$-set $K :=A \backslash H$
with $A$ subgroup of $H$, and a cochain $\mu \in C^2(H, \Map(K, \Cx))$ such that $\delta_H \mu = \pi^* \eta$ 
with $\pi^* \eta(k; h_1,h_2,h_3)= \eta(h_1,h_2,h_3)$ (see \cite[\S 3.3]{Uribe}).
The skeletal tensor category of the tensor category $\CC^*_\MM= Fun_\CC(\MM,\MM)$ is equivalent to one of 
the form $\VV(\widehat{H},\widehat{\eta})$ whenever $A$ is normal and abelian in $H$ \cite{Naidu} and if there exists a cochain $\gamma \in C^1(H, \Map(K, \Cx))$ such that $\delta_H \gamma= \delta_K \mu$. In particular this implies that the cohomology class of $\eta$
belongs to the subgroup of $H^3(H, \Cx)$ defined by
$$\Omega(H;A) := \ker \left( \ker \left( H^3(H, \Cx) \to E^{0,3}_\infty \right) \to E^{1,2}_\infty  \right),$$
which fits into the short exact sequence \cite[Cor. 3.2]{Uribe}
 $$0 \to E^{3,0}_\infty \to \Omega(G;A) \to E^{2,1}_\infty \to 0$$
where $E_n^{*,*}$ denotes the $n$-th page of the Lyndon-Hochschild-Serre spectral sequence associated to the group extension
$1 \to  A \to H \to K \to 1$.

Let us bring the last piece of notation. Denote the dual group ${{\mathbb{A}}} := \Hom(A, \complex^*)$ and consider 
cocycles $F \in Z^2(K, A)$ and $\widehat{F} \in Z^2(K, {{\mathbb{A}}})$. Denote by $G= A \rtimes_F K$ and
$\widehat{G}= K \ltimes_{\widehat{F}} {{\mathbb{A}}}$ the groups defined by the multiplication laws 
$$(a_1,k_1) (a_2,k_2) := (a_1 ({}^{k_1}a_2) F(k_1,k_2),k_1k_2)$$
$$(k_1, \rho_1) \cdot (k_2, \rho_2) := (k_1k_2, (\rho_1^{k_2}) \rho_2 \widehat{F}(k_1,k_2))$$
respectively.  The necessary and sufficient  conditions for two pointed fusion categories to be Morita equivalent are the following (cf. \cite{Naidu}):

\begin{theorem}\cite[Thm. 5.9]{Uribe}  \label{classification theorem}
 Let $H$ and $\widehat{H}$ be finite groups, $\eta \in Z^3(H, \complex^*)$ and $\widehat{\eta} \in Z^3(\widehat{H}, \complex^*)$. Then the tensor categories $Vect(H,\eta)$ and $Vect(\widehat{H}, \widehat{\eta})$ are weakly Morita equivalent if and only if the following
  conditions are satisfied:
 \begin{itemize}
 \item There exist isomorphisms of groups 
 $$\phi : G= A \rtimes_F K \stackrel{\cong}{\to} H \ \ \ \ 
 \widehat{\phi} : \widehat{G}= K \ltimes_{\widehat{F}} {{\mathbb{A}}} \stackrel{\cong}{\to} \widehat{H}$$
 for some finite group $K$ acting on the abelian group $A$,
 with cocycles $F \in Z^2(K, A)$ and $\widehat{F} \in Z^2(K, {{\mathbb{A}}})$. 
 \item There exists $\epsilon : K^3 \to \Cx$ such that $\widehat{F} \wedge F =\delta_K \epsilon$.
 \item The cohomology classes satisfy the equations $[\omega]=[ \phi^* \eta ]$ and 
 $[\widehat{\omega}]=[\widehat{\phi}^*\widehat{\eta}]$ with 
 \begin{align*}
   \omega((a_1,k_1),(a_2,k_2),(a_3,k_3)) := & \widehat{F}(k_1,k_2)(a_3) \ \epsilon(k_1,k_2,k_3)\\
    \widehat{\omega}((k_1, \rho_1),  (k_2,\rho_2),(k_3 ,\rho_3)) :=& \epsilon(k_1,k_2,k_3) \ \rho_1(F(k_2,k_3)).
   \end{align*}
\end{itemize}
\end{theorem}

Note that ${{\mathbb{A}}}$ and $A$ are (non-canonically) isomorphic as $K$-modules, and therefore
both $H$ and $\widehat{H}$ could be seen as extensions of $K$ by $A$. In order to calculate
all possible Morita equivalences, we will analize the Morita equivalences that appear while fixing the 
group $K$ and the $K$-module $A$.

Let us recall the equivalence classes of normal an abelian subgroups of the groups of order 8. We will
say that two subgroups are equivalent if there is an automorphism of the group that maps one to the other.
The following table contains the information on these subgroups:
\begin{center}
 \begin{tabular}{|c|c|c|c|c|c|}
\hline
Isomorphic to&$D_8$ & $\integer/{4}\oplus\integer/{2}$ & $Q_{8}$ & $\integer/{8}$\tabularnewline
\hline
 $\integer/2$ &$\left\langle a^{2}\right\rangle $ & $\left\langle (2,0)\right\rangle $ & $\left\langle -1\right\rangle $ & $\left\langle 4\right\rangle $\tabularnewline
\hline
 $\integer/2$ & & $\left\langle (0,1)\right\rangle, $
$\left\langle (2,1)\right\rangle $ & &
\tabularnewline
\hline
$\integer/2\times\integer/2$ &$\left\langle a^{2},ba^{3}\right\rangle, \left\langle a^{2},b\right\rangle $ &$\left\langle (2,0)\right\rangle \oplus\left\langle (0,1)\right\rangle $ & &
\tabularnewline
\hline
\integer/4&$\left\langle a\right\rangle $ & $\left\langle (1,0)\right\rangle,   \left\langle (1,1)\right\rangle$ &$\left\langle i\right\rangle $,$\left\langle j\right\rangle $,$\left\langle k\right\rangle $ &$\left\langle 2\right\rangle $,$\left\langle 6\right\rangle $ 
\tabularnewline
\hline
\end{tabular}
\end{center}
In the group $(\integer/2)^3$ all subgroups of order 2 are isomorphic to $\langle (0,0,1) \rangle$ and all subgroups of order 4 are isomorphic
to $\langle (0,1,0), (0,0,1) \rangle$.

Now let us outline the procedure that we will follow. We will fix the groups $K$ and $A$, we will take the groups
that are extensions of $K$ by $A$ and we will take explicit choices of subgroups
from the table above that provide the extensions. Then we will calculate the relevant terms of the second page of the Lyndon-Hochschild-Serre
spectral sequence, which are the same for all extensions of $K$ by $A$, and we will calculate the third page for each extension $0 \to A \to H \to K\to 1$. Then we will determine the cohomology class
of the 2-cocycle $F$ that makes $H \cong A \rtimes_F K$ and we will calculate the cohomology classes in $\Omega(H;A)$.
With this information and Theorem \ref{classification theorem} we will determine the Morita equivalence classes of pointed fusion categories for groups 
that are extensions of $K$ by $A$.

\subsection{$K= \integer/2$ and $A=\integer/4$ with trivial action}

The two possible extensions are $\integer/8$ and $\integer/4 \times \integer/2$ with the following choices of subgroups:
\begin{align*} 
1 \longrightarrow\langle (1,0)\rangle \longrightarrow \integer/4 & \times \integer/2 \longrightarrow \integer/2 \longrightarrow 1  \\
1 \longrightarrow\langle 2\rangle \longrightarrow & \integer/8  \longrightarrow \integer/2  \longrightarrow 1.
\end{align*}
For the group $\integer/4 \times \integer/2$ the relevant terms of the second page of the LHS spectral sequence are:
\renewcommand{\arraystretch}{2}
\begin{align*}\label{rob}
\newcommand*{\tempb}{\multicolumn{1}{|c}{}}
\begin{array}{cccccccccc}
3 & \tempb & \mathbb{Z}/4=\langle u^2\rangle\\
2 & \tempb & 0& 0&0\\
1 & \tempb & \Z/4& \Z/2 & \mathbb{Z}/2 =\langle uv\rangle\\
0 & \tempb & \mathbb{C}^*& \mathbb{Z}/2 & 0 & \mathbb{Z}/2 =\langle v^2\rangle& 0 \\ \cline{2-7}
& & 0~~~& 1~~~ & 2~~~ & 3~~~ & 4
\end{array}
\end{align*}
Since $H^2(\integer/2, \integer/4) = \integer /2$ we have that 
$\integer/8 \cong \integer/2 \ltimes_{uv} {\integer/4}$, and therefore the relevant terms of the third page of the LHS spectral sequence
for $\integer/8$ are:

 \begin{tikzpicture}
\matrix [matrix of math nodes,row sep=6mm]
{
 3 &  [5mm]  |(a)|  \integer/4=\langle s^2\rangle/( 4s^2) & [5mm]   & [5mm]  & [5mm] & [5mm] & [5mm] \\\
2 & |(b)| 0 & |(c)|  0  &  & & & \\
1&  \integer/4 & |(d)|  \integer/2& |(e)|\integer/2 =\langle 4s^2\rangle &  & & \\
0& \Cx &  \integer/2 & |(f)| 0 & |(g)|  \integer/2& 0 \\
& 0 & 1 & 2 & 3& 4&\\
};
\tikzstyle{every node}=[midway,auto,font=\scriptsize]
\draw[thick] (-4.5,-1.7) -- (-4.5,2.8) ;
\draw[thick] (-4.5,-1.7) -- (5.0,-1.7) ;
%\draw[-stealth] (b) -- node {$\delta_G$} (c);
\draw[-stealth] (d) -- node {$\cong$} (g);
%\draw[-stealth] (d) -- node {$\delta_G$} (e);
%\draw[-stealth] (f) -- node {$\delta_K$} (e);
%\draw[-stealth] (f) -- node {$\delta_G$} (g);
\end{tikzpicture}

We get the split short exact sequences 
\begin{align*}0 \to  \langle v^2\rangle  &\to \Omega(\integer/4 \times \integer/2; \langle (1,0)\rangle) \to \langle uv\rangle \to 0\\
& 0 \to \Omega(\integer/8; \langle 4 \rangle) \stackrel{\cong}{\to} \langle 4s^2 \rangle \to 0
 \end{align*} and we conclude that the only Morita equivalences that appear 
are:
\begin{align*}
Vect(\integer/8,0)\simeq_M Vect(\integer/4\times\integer/2,uv)\simeq_MVect(\integer/4\times\integer/2,uv+v^2).
\end{align*}

\subsection{$K=\integer/2$ and $A=\integer/4$ with non-trivial action.}
In this case the possible extensions are:
\begin{align*}
1  \longrightarrow \langle a\rangle\longrightarrow D_8 \longrightarrow \integer/2  \longrightarrow 1\\
1  \longrightarrow \langle i\rangle\longrightarrow Q_8\longrightarrow \integer/2  \longrightarrow 1.
\end{align*}
 where $D_8$ is the trivial extension of $\integer/2$ by the non-trivial module $\integer/4$. For $D_8$ the relevant elements of the second page of the LHS spectral sequence are:
    
    \begin{tikzpicture}
\matrix [matrix of math nodes,row sep=6mm]
{
 3 &  [5mm]  |(a)|  \integer/4= \langle \gamma \rangle & [5mm]   & [5mm]  & [5mm] & [5mm] & [5mm] \\\
2 & |(b)| 0 & |(c)|  0  &  & & & \\
1&  \integer/2 & |(d)|  \integer/2& |(e)| \integer/2= \langle \beta \rangle &  & & \\
0& \Cx &  \integer/2 & |(f)| 0 & |(g)|  \integer/2=\langle \alpha \rangle & 0 \\
& 0 & 1 & 2 & 3& 4&\\
};
\tikzstyle{every node}=[midway,auto,font=\scriptsize]
\draw[thick] (-4.5,-1.7) -- (-4.5,2.8) ;
\draw[thick] (-4.5,-1.7) -- (5.0,-1.7) ;
%\draw[-stealth] (b) -- node {$\delta_G$} (c);
%\draw[-stealth] (d) -- node {$(\delta_K)^{-1}$} (c);
%\draw[-stealth] (d) -- node {$\delta_G$} (e);
%\draw[-stealth] (f) -- node {$\delta_K$} (e);
%\draw[-stealth] (f) -- node {$\delta_G$} (g);
\end{tikzpicture}

\noindent where the generators of the cohomology of $D_8$ were defined in \eqref{generators D8}. Note that the cohomology
class $\alpha$ is the pullback of the generator of $H^4(\langle b \rangle, \integer)$ and $\beta$ is the cohomology
class in $H^2(\integer/2, \integer/4)$ that classifies the non-trivial extension determined by $Q_8$.

   The second page of the LHS spectral sequence for the extension  $Q_8 \cong \integer/2 \ltimes_\beta {\integer/4}$ becomes:
   
    \begin{tikzpicture}
\matrix [matrix of math nodes,row sep=6mm]
{
 3 &  [5mm]  |(a)|  \integer/4= \langle t \rangle/(4t) & [5mm]   & [5mm]  & [5mm] & [5mm] & [5mm] \\\
2 & |(b)| 0 & |(c)|  0  &  & & & \\
1&  \integer/2 & |(d)|  \integer/2& |(e)| \integer/2= \langle 4t\rangle/(8t) &  & & \\
0& \Cx &  \integer/2 & |(f)| 0 & |(g)|  \integer/2& 0 \\
& 0 & 1 & 2 & 3& 4&\\
};
\tikzstyle{every node}=[midway,auto,font=\scriptsize]
\draw[thick] (-4.5,-1.7) -- (-4.5,2.8) ;
\draw[thick] (-4.5,-1.7) -- (5.0,-1.7) ;
%\draw[-stealth] (b) -- node {$\delta_G$} (c);
\draw[-stealth] (d) -- node {$\cong$} (g);
%\draw[-stealth] (d) -- node {$\delta_G$} (e);
%\draw[-stealth] (f) -- node {$\delta_K$} (e);
%\draw[-stealth] (f) -- node {$\delta_G$} (g);
\end{tikzpicture}

\noindent   where the second differential $d_2: E_2^{1,1} \stackrel{\cong}{\to} E_2^{3,0}$ is an isomorphism since $H^3(Q_8,\Cx)=\langle t \rangle/(8t)$. 

We get the split short exact sequences 
\begin{align*}0 \to  &\langle \alpha \rangle  \to \Omega(D_8; \langle a \rangle) \to \langle \beta \rangle \to 0\\
 &0 \to \Omega(Q_8; \langle i \rangle) \stackrel{\cong}{\to} \langle 4t  \rangle \to 0
 \end{align*} and we conclude that the only Morita equivalences that appear 
are:
   \begin{align*}
   Vect(Q_8, 0) \simeq_M Vect(D_8,\beta) \simeq_M  Vect(D_8,\beta \oplus \alpha). \end{align*} 

\subsection{$K=\integer/4$ and $A=\integer/2$.}
The two extensions are:
\begin{align*}
1  \longrightarrow\langle (0,1)& \rangle \longrightarrow \integer/4 \times \integer/2 \to \integer/4  \longrightarrow 1\\
1  \longrightarrow&\langle 4 \rangle \longrightarrow \integer/8 \longrightarrow \integer/4  \longrightarrow 1,
\end{align*}
and the relevant terms of the second page of the LHS spectral sequence for $\integer/4 \times \integer/2$ are:
\renewcommand{\arraystretch}{2}
$$\newcommand*{\tempb}{\multicolumn{1}{|c}{}}
\begin{array}{cccccccccc}
3 & \tempb & \mathbb{Z}/2=\langle v^2\rangle\\
2 & \tempb & 0& 0&0\\
1 & \tempb &\mathbb{Z}/2& \mathbb{Z}/2& \mathbb{Z}/2=\langle uv\rangle\\
0 & \tempb & \mathbb{C}^*& \mathbb{Z}/4 & 0 & \mathbb{Z}/4=\langle u^2\rangle & 0 \\ \cline{2-7}
 & & 0& ~~~~1~~~~ & ~~~~2~~~~ & ~~~~3~~~~ & 4
\end{array}
$$
We have again that  $\integer/8=\integer/4 \rtimes_{uv}\integer/2$, and the relevant terms of the third page of the LHS spectral sequence
associated to the extension of $\integer/8$ are:
\renewcommand{\arraystretch}{2}
$$\newcommand*{\tempb}{\multicolumn{1}{|c}{}}
\begin{array}{cccccccccc}
3 & \tempb & \mathbb{Z}/2=\langle s^2\rangle/(2s^2)\\\
2 & \tempb & 0& 0&0\\
1 & \tempb &\mathbb{Z}/2& 0 & \mathbb{Z}/2=\langle2s^2\rangle/(4s^2)\\
0 & \tempb & \mathbb{C}^*& \mathbb{Z}/4 & 0 & \mathbb{Z}/2=\langle 4s^2\rangle & 0 \\ \cline{2-7}
 & & 0& ~~~~1~~~~ & ~~~~2~~~~ & ~~~~3~~~~ & 4
\end{array}
$$
We obtain the short exact sequences 
\begin{align*}
&0 \to  \langle v^2 \rangle  \to \Omega(\integer/4 \times \integer/2; \langle (0,1) \rangle) \to \langle uv \rangle \to 0\\
 0 &\to \langle 4s^2 \rangle /(8s^2)\to \Omega(\integer/8; \langle 2 \rangle) \to \langle 2s^2 \rangle /(4s^2) \to 0
 \end{align*} 
 where $\Omega(\integer/8; \langle 2 \rangle) \cong \integer/4$. We conclude that the only Morita equivalences that appear 
are:

\begin{align*}
&Vect(\integer/8,0)\simeq_M Vect(\integer/4\times \integer/2,uv)\\
Vect(\integer/8,4s^2)\simeq_M & Vect(\integer/4\times\integer/4,uv+u^2) \simeq_M Vect(\integer/4\times\integer/4,uv+3u^2)
\end{align*} 

\subsection{$K=\integer/2$ and $A= \integer/2 \times \integer/2$ with trivial action}
The relevant extensions are:
\begin{align*}
1  \longrightarrow\langle (0,0,1), (0,1,0) \rangle \longrightarrow (\integer/2)^3 \to \integer/2  \longrightarrow 1\\
1  \longrightarrow\langle (2,0),(0,1) \rangle \longrightarrow \integer/4 \times \integer/2 \longrightarrow \integer/2  \longrightarrow 1.
\end{align*}
The relevant terms for the second page $E_2^{p,q}=H^p(\integer/2,H^q(\integer/2 \times \integer/2, \complex^*)$ of the LHS spectral sequence are:
\renewcommand{\arraystretch}{2}
\begin{align*}
\newcommand*{\tempb}{\multicolumn{1}{|c}{}}
\begin{array}{ccccccccccc}
3 & \tempb & (\mathbb{Z}/2)^3\\
2 & \tempb & \langle y^2z+yz^2 \rangle&&\langle xyz \rangle\\
1 & \tempb & \langle y^2,z^2 \rangle &&\langle xy,xz \rangle&& \langle x^2y,x^2z \rangle\\
0 & \tempb & \mathbb{C}^*&& \langle x^2 \rangle && 0 && \langle x^4 \rangle&& 0 \\ \cline{2-11}
& & 0&& 1 && 2 && 3 && 4
\end{array}
\end{align*}
where on the base we have $H^*(\integer/2, \field_2) \cong \field_2[x]$ and on the fiber
$H^*(\integer/2 \times \integer/2, \field_2) \cong \field_2[y,z] $. The cohomology classes
 that appear on the 0-th row an the 0-th column are the ones that are annihilated  by the operation $Sq^1$.

The second differential of the LHS spectral sequence of the extension of $\integer/4 \times \integer/2$,
whose $k$-invariant is the class $x^2$, maps \begin{align*}
 y^2z + yz^2 \mapsto x^2y, \ \ xz \mapsto x^4, \ \ xy\mapsto 0,\ \ xyz\mapsto x^3y.
 \end{align*} Therefore the relevant terms of the third page of the LHS spectral sequence for $\integer/4 \times \integer/2$ are
 (using the base defined in \S \ref{subsection Z4xZ2}):
\renewcommand{\arraystretch}{2}
\begin{align*}
\newcommand*{\tempb}{\multicolumn{1}{|c}{}}
\begin{array}{ccccccccccc}
3 & \tempb & \langle u^2, uv, v^2 \rangle /(2u^2)\\
2 & \tempb & 0&&0\\
1 & \tempb & \langle u,v \rangle /(2u) &&\integer/2&& \langle 2u^2 \rangle\\
0 & \tempb & \mathbb{C}^*&& \langle 2u \rangle && 0 &&  0&& 0 \\ \cline{2-11}
& & 0&& 1 && 2 && 3 && 4
\end{array}
\end{align*}
We get the split short exact sequences 
\begin{align*}0 \to  &\langle x^4 \rangle  \to \Omega((\integer/2)^3; \langle (0,0,1),(0,1,0) \rangle) \to \langle x^2z^2, x^2y^2 \rangle \to 0\\
 &0 \to  \Omega(\integer/4 \times \integer/2; \langle (2,0),(0,1) \rangle) \stackrel{\cong}{\to} \langle 2u^2  \rangle \to 0
 \end{align*} where all the non-trivial classes in $\langle x^2z^2, x^2y^2 \rangle$ define an extension isomorphic to $\integer/4 \times \integer/2$,
 and $\integer /4 \times \integer/2 \cong \integer/2 \ltimes_{x^2y^2} (\integer/2 \times \integer/2)$.
We conclude that the only Morita equivalence that appear is:
\begin{align*}
Vect(\integer/4 \times\integer/2 ,0)\simeq_M Vect((\integer/2)^3,x^2z^2).
\end{align*} 

 \subsection{$K=\integer/2$ and $A=\integer/2 \times \integer/2$ with non-trivial action}

Let us take the action of $\integer/2$ that flips the coordinates in $\integer/2 \times \integer/2$, therefore the only extension
is equivalent to the following one:
$$0 \longrightarrow \langle b, ba^2 \rangle \longrightarrow D_8 \longrightarrow \integer/2 \longrightarrow 1.$$
Since there is up to isomorphism only one extension, this choice of $K$ and $A$ does not produce any Morita equivalence
between non-equivalent pointed fusion categories.

\subsection{$K=\integer/2\times \integer/2$ and $A=\integer/2$}
In this last case the relevant extensions are:
\begin{align*}
1 \longrightarrow \langle (0,0,1) \rangle  \longrightarrow& (\integer/2)^3 \longrightarrow \integer/2 \times \integer/2 \longrightarrow 1\\
1 \longrightarrow \langle (2,0) \rangle \longrightarrow \integer/4 &\times \integer/2 \longrightarrow  \integer/2 \times \integer/2 \longrightarrow 1\\
1 \longrightarrow \langle a^2 \rangle \longrightarrow & D_8 \longrightarrow\integer/2 \times \integer/2 \longrightarrow 1\\
1 \longrightarrow \langle -1 \rangle \longrightarrow &Q_8 \longrightarrow \integer/2 \times \integer/2 \longrightarrow 1.
\end{align*}
The relevant terms of the second page $E_2^{p,q}=H^p(\integer/2 \times \integer/2, H^q(\integer/2, \complex^*))$
of the LHS spectral sequence are:
\renewcommand{\arraystretch}{2}
$$\newcommand*{\tempb}{\multicolumn{1}{|c}{}}
\begin{array}{cccccccccc}
3 & \tempb & \langle z^4\rangle\\
2 & \tempb & 0& 0 & 0 \\
1 & \tempb &\langle z^2 \rangle& \langle zx, zy\rangle&\langle zx^2,zy^2,zxy\rangle\\
0 & \tempb & \mathbb{C}^*& \langle x^2, y^2\rangle & \langle x^2y+xy^2\rangle&\langle x^4,y^4, x^2y^2\rangle& \langle  x^4y+x^3y^2, x^2y^3+xy^4 \rangle  \\ \cline{2-7}
 & & 0& ~~~~1~~~~ & ~~~~2~~~~ & ~~~~3~~~~ & 4
\end{array}
$$
where we have assumed that $H^*(\integer/2 \times \integer/2, \field_2)=\field_2[x,y]$ for the base, $H^*(\integer/2 , \field_2)=\field_2[z]$
for the fiber, and the elements on the 0-th row and the 0-th column are the classes that get annihilated by the operation $Sq^1$.

The calculation of the second differential between the first row and the 0-th row 
$$d_2^G : H^p((\integer/2)^2; \field_2) \otimes H^1(\integer/2; \field_2) \to H^{p+2}((\integer/2)^2, \complex^*)$$
depends on the $k$-invariant of the extension $G$. In \S \ref{subsection dihedral} we explain the case for $D_8$ on which
the differential was $$d_2^{D_8}(zp(x,y)) = Sq^1((xy+x^2)p(x,y))$$ and following the same argument we get that the second 
differentials for the other two groups are
\begin{align*}
&d_2^{\integer/4 \times \integer/2}(zp(x,y)) = Sq^1((x^2)p(x,y))\\
&d_2^{Q_8}(zp(x,y)) = Sq^1((x^2+xy+x^2)p(x,y)).
\end{align*}

\subsubsection{} For the group $(\integer/2)^3$ we obtain a split extension
$$0 \to \langle x^4,y^4, x^2y^2\rangle \to \Omega((\integer/2)^3; \langle (0,0,1) \rangle) \to \langle zx^2,zy^2,zxy\rangle \to 0$$
where the classes $x^2,x^2+y^2,y^2$ in $H^2((\integer/2)^2, \field_2)$ define an extension isomorphic to $\integer/4 \times \integer/2$,
the classes $xy, xy+x^2,xy+y^2$ define an extension isomorphic to $D_8$ and the class $x^2+xy+y^2$ define an extension isomorphic to $Q_8$.
Therefore it is enough to analize the cases determined by the classes $x^2$, $xy+x^2$ and $x^2+xy+y^2$. 

\subsubsection{} For the group $\integer/4 \times \integer/2$ we have that it is isomorphic to the group
 $\integer/2\ltimes_{x^2}(\integer/2\times\integer/2)$ and the relevant terms of the third page of the LHS are:

\renewcommand{\arraystretch}{2}
\begin{align*}
\newcommand*{\tempb}{\multicolumn{1}{|c}{}}
\begin{array}{cccccccccc}
3 & \tempb & \langle z^4\rangle\\
2 & \tempb & 0& 0\\
1 & \tempb & \langle z^2 \rangle & 0& \langle x^2z,y^2z\rangle \\
0 & \tempb & \mathbb{C}^*& \langle x^2,y^2\rangle & \langle x^2y+xy^2 \rangle & \dfrac{\langle x^4,x^2y^2,y^4\rangle}{\langle x^2y^2,x^4\rangle}& \integer/2 \\ \cline{2-7}
& & 0~~~& 1~~~ & 2~~~ & 3~~~ & 4
\end{array}
\end{align*}
Following the notations of \S \ref{subsection Z4xZ2} we know that $E_3^{3,0}$ is generated by $v^2$ and
$E_3^{2,1}$ is generated by $uv$ and $2u^2$, where $uv$ corresponds to the class $y^2z$ and $2u^2$ corresponds to the class $x^2z$.
Therefore we obtain
$$0 \to \langle v^2 \rangle \to \Omega(\integer/4 \times \integer/2; \langle (2,0) \rangle ) \to \langle uv, 2u^2 \rangle \to 0$$
and since all the classes in $E_3^{2,1}$ induce extensions isomorphic to $\integer/4 \times \integer/2$, we conclude
that the only Morita equivalences that we obtain in this case are:
\begin{align*}
Vect(\integer/4 \times \integer/2, 0) &\simeq_M Vect((\integer/2)^3,x^2z^2)\\
Vect(\integer/4 \times \integer/2, v^2) & \simeq_M Vect((\integer/2)^3,x^2z^2+y^4).
\end{align*}

\subsubsection{} The calculation of the relevant terms of the third page of the LHS spectral sequence for the group 
$D_8 \cong \integer/2 \ltimes_{xy+x^2}(\integer/2 \times \integer/2 )$ was done in \S \ref{subsection dihedral}
and therefore we obtain the short exact sequence
$$0 \to \langle \alpha, \beta \rangle \to \Omega(D_8; \langle a^2 \rangle) \to \langle 2\gamma \rangle \to 0.$$
Since $E_3^{2,1}$ is generated by the class $z(xy+y^2)$ and $xy+x^2$ is the $k$-invariant for $D_8$, we only get
the following Morita equivalences:
\begin{align*}
Vect(D_8, &0) \simeq_M Vect((\integer/2)^3,x^2yz+xy^2z+xyz^2+ x^2z^2)\\
Vect(D_8, \alpha&) \simeq_M Vect((\integer/2)^3,y^4 + x^2yz+xy^2z+xyz^2+ x^2z^2)\\
Vect(D_8, \beta&) \simeq_M Vect((\integer/2)^3,x^4+ x^2yz+xy^2z+xyz^2+ x^2z^2).
\end{align*}

\subsubsection{} The relevant terms of the third page for the LHS spectral sequence for the group
$Q_8 \cong \integer/2 \ltimes_{x^2+xy+y^2}(\integer/2 \times \integer/2 )$ are:
\renewcommand{\arraystretch}{2}
\begin{align*}
\newcommand*{\tempb}{\multicolumn{1}{|c}{}}
\begin{array}{cccccccccc}
3 & \tempb & \langle z^4\rangle\\
2 & \tempb & 0& 0\\
1 & \tempb & 0 & 0 &\langle z(x^2+ xy+y^2)\rangle \\
0 & \tempb & \mathbb{C}^*&\langle x^2,y^2\rangle & 0 &\dfrac{\langle x^4,y^4, x^2y^2\rangle}{\langle x^4,y^4\rangle}& 0 \\ \cline{2-7}
& & 0~~~& 1~~~ & 2~~~ & 3~~~ & 4
\end{array}
\end{align*}
Since $H^3(Q_8,\complex^*)\cong \integer\langle  t \rangle /(8t)$ we know that $E_3^{0,3}\cong \integer\langle  t \rangle /(2t)$,
$E_3^{2,1}\cong \integer\langle  2t \rangle /(4t)$ and $E_3^{3,0}\cong \integer\langle  4t \rangle /(8t)$. Therefore
we get that $\Omega(Q_8;\langle -1 \rangle) = \integer \langle 2t \rangle /(8t)$ which fits into the short exact sequence:
$$0 \to \integer\langle  4t \rangle /(8t) \to \Omega(Q_8;\langle -1 \rangle) \to \integer\langle  2t \rangle /(4t) \to 0.$$
Since the $E_3^{2,1}$ is generated by the class $z(x^2+xy+y^2)$ and $x^2+xy+y^2$ is the $k$-invariant of $Q_8$, we only get 
the following Morita equivalences:
\begin{align*}
Vect(Q_8, 0)& \simeq_M Vect((\integer/2)^3,x^2yz+xy^2z+xyz^2+ x^2z^2+y^2z^2)\\
Vect(Q_8, 4t) \simeq_M & Vect((\integer/2)^3,x^2yz+xy^2z+xyz^2+ x^2z^2+y^2z^2+x^2y^2).\\
\end{align*}

\subsection{} We conclude that the only weak Morita equivalences between pointed fusion categories
of global dimension 8 are the ones that appear on each row of the following table:
\begin{center}
\begin{tabular}{ ||c|c|c|c|c|| } 
 \hline 
 $(\integer/2)^3$ & $\integer/4 \times \integer/2$ & $\integer/8$ &  $D_8$ & $Q_8$\\[0.5ex] 
 \hline\hline
 $orb(x^2y^2)$ & $\{0\}$ &&& \\ \hline
  $orb(x^4+y^2z^2)$ & $\{v^2\}$ &&& \\ \hline
   $orb(x^2yz+$$xy^2z+xyz^2)$ &  &  &$\{0\}$& \\ \hline
 $orb(x^4+ x^2yz+xy^2z+xyz^2)$ &  &  &$\{\alpha+\beta, \beta \}$& $\{0\}$\\ \hline
 $orb(x^4+ x^2yz+xy^2z+xyz^2+y^2z^2)$ &  &  &$\{\alpha \}$& \\ \hline
  $orb(x^2yz+xy^2z+xyz^2+x^2y^2+x^2z^2+y^2z^2)$ &&&& $\{4t\}$ \\ \hline
  & $orb(uv)$ & $\{0\}$ &  & \\ \hline
  & $orb(uv + u^2)$ & $\{4s^2\}.$ && \\ \hline
\end{tabular}
\end{center}

Therefore there are only 47-9=36 Morita equivalence classes of pointed fusion categories of global dimension 8.

\subsection{Twisted Drinfeld double}

  The center of a fusion category is again a fusion category and it is moreover braided.
   In \cite[Thm. 3.1]{W-GT-fusion-cat} it is shown that two tensor categories are weakly Morita
      equivalent if and only if their centers are braided equivalent.
 In particular, if $\MM$ is an indecomposable module
 category over $\CC$, there is a canonical equivalence of braided tensor categories
 $\ZZ(\CC) \simeq \ZZ(\CC_\MM^*)$ \cite[Prop. 2.2]{Ost-2}.  The center $\ZZ(Vect(H, \eta))$ of the tensor category
  $Vect(H, \eta)$ contains the information necessary for constructing the braided quasi-Hopf algebra that is known as the twisted
  Drinfeld double $D^\eta(H)$ of the group $H$ twisted by $\eta$ (see \cite[\S 3.2]{Dijkgraaf}), and moreover
 the centers of two pointed fusion categories are braided equivalent if and only if the associated twisted Drinfeld doubles
 are isomorphic as braided quasi-Hopf algebras.
 
Therefore the twisted Drinfeld doubles $D^\eta(H)$ and $D^{\widehat{\eta}}(\widehat{H})$
are isomorphic as braided quasi-Hopf algebras if and only if the pointed fusion categories 
 $Vect(H, \eta)$ and  $Vect(\widehat{H}, \widehat{\eta})$ are weakly Morita equivalent.
 
 We conclude that there are 38 isomorphism classes of twisted Drinfeld doubles of 
 groups of order 8, 18 of them corresponding to twisted Drinfeld doubles which are commutative as algebras
 and 20 corresponding  to twisted Drinfeld doubles which are non-commutative as algebras.
 The commutative twisted Drinfeld doubles are the ones constructed from the groups $\integer/8$, $\integer/4 \times \integer/2$ with any associator and from the group $(\integer/2)^3$ with associator any cohomology class
 not containing the element $Sq^1(xyz)= x^2yz + xy^2z + xyz^2$. The non-commutative twisted
 Drinfeld doubles are the ones constructed from the groups $D_8$, $Q_8$ with any associator 
 and from the group $(\integer/2)^3$ with associator any cohomology class
 containing the element $Sq^1(xyz)= x^2yz + xy^2z + xyz^2$. 
 
   The isomorphism classes of twisted Drinfeld doubles of groups of order 8 was
   carried out with complete different methods in \cite{Mason} whenever the algebra structure
   was commutative and in \cite{Goff} whenever the algebra structure was not commutative.
   Our results are compatible with the ones of the previous two references and independent of them.

   % \bibliographystyle{alpha}
%\bibliography{Categorical}
  % \end{document}
   
   %\bibliographystyle{alpha} 
\bibliographystyle{abbrv} 
\def\cprime{$'$} \def\cprime{$'$}

\end{document}